\documentclass[11pt,reqno]{amsart}

\usepackage[margin=1in]{geometry}

\usepackage{amsmath,amssymb,amscd,amsxtra}
\usepackage[dvips]{graphics}
\usepackage{float}

\usepackage{graphicx}
\usepackage{psfrag}

\usepackage[format=hang,indention=0pt,figurename=Fig.]{caption}

\usepackage{amsfonts}
\usepackage{amsthm}
\usepackage{enumerate}
\usepackage{bbm}

\usepackage{hyperref}

\newcommand{\nD}{\mathbb D}
\newcommand{\nE}{\mathbb E}
\newcommand{\nF}{\mathbb F}
\newcommand{\nG}{\mathbb G}

\newcommand{\nR}{\mathbb R}

\newcommand{\Vb}[1]{\mathbf{{#1}}}

\newtheorem*{theorem*}{Theorem}

\theoremstyle{definition}

\begin{document}

\title{Outer median triangles}

%%%  all figures made in v8 Outer_Median_Triangle.nb

%----------Author 1
\author[\'A. B\'enyi]{\'Arp\'ad B\'enyi}

\address{%
Department of Mathematics \\
Western Washington University \\
516 High Street, Bellingham \\
Washington 98225, USA}

\email{Arpad.Benyi@wwu.edu}

\thanks{This work is partially supported by a grant from the Simons Foundation (No.~246024 to \'Arp\'ad B\'enyi).}

%----------Author 2
\author[B. \'{C}urgus]{Branko \'{C}urgus}
\address{%
Department of Mathematics \\
Western Washington University \\
516 High Street, Bellingham \\
Washington 98225, USA}

\email{Branko.Curgus@wwu.edu}



\date{\today}

%%% ----------------------------------------------------------------------

\begin{abstract}
We define the notions of outer medians and outer median triangles. We show that outer median triangles enjoy similar properties to that of the median triangle.
\end{abstract}

\maketitle

A \emph{median of a triangle} is a line segment that connects a vertex of the triangle to the midpoint of the opposite side. The three medians of a given triangle interact nicely with each other to yield the following properties:
\begin{enumerate}[{\rm (a)}]
\item  \label{Pra}
The medians \emph{intersect in a point} interior to the triangle,  called the \emph{centroid}, which divides each of the medians in the ratio 2:1.
\item \label{Prb}
The medians form a new triangle, called the \emph{median triangle}.
\item \label{Prc}
The area of the median triangle is 3/4 of the area of the given triangle in which the medians were constructed.
\item \label{Prd}
The median triangle of the median triangle is similar to the given triangle with the ratio of similarity 3/4.
\end{enumerate}

\begin{figure}[H]

\setlength{\abovecaptionskip}{5pt}%
\setlength{\belowcaptionskip}{-5pt}%

\psfrag{A}[][]{\begin{picture}(0,0)
            \put(-3,1){\makebox(0,0)[l]{$A$}}
                        \end{picture}}

\psfrag{B}[][]{\begin{picture}(0,0)
            \put(-6,-3){\makebox(0,0)[l]{$B$}}
                        \end{picture}}
\psfrag{C}[][]{\begin{picture}(0,0)
            \put(-2,-3){\makebox(0,0)[l]{$C$}}
                        \end{picture}}

\psfrag{At}[][]{\begin{picture}(0,0)
            \put(-5,-5){\makebox(0,0)[l]{$A_{1/2}$}}
                        \end{picture}}

\psfrag{Bt}[][]{\begin{picture}(0,0)
            \put(0,1){\makebox(0,0)[l]{$B_{1/2}$}}
                        \end{picture}}

\psfrag{Ct}[][]{\begin{picture}(0,0)
            \put(-16,2){\makebox(0,0)[l]{$C_{1/2}$}}
                        \end{picture}}

\resizebox{0.7\linewidth}{!}{%
 \includegraphics{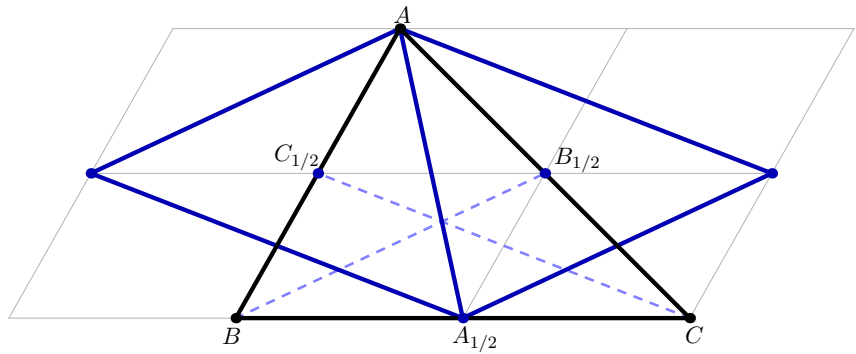}}
    \caption{
A ``proof'' of Properties (\ref{Prb}) and (\ref{Prc})
        }
\label{famt}

\end{figure}

\begin{figure}[H]

\setlength{\abovecaptionskip}{5pt}%
\setlength{\belowcaptionskip}{-5pt}%

\psfrag{A}[][]{\begin{picture}(0,0)
            \put(-4,1){\makebox(0,0)[l]{$A$}}
                        \end{picture}}

\psfrag{B}[][]{\begin{picture}(0,0)
            \put(-6,-3){\makebox(0,0)[l]{$B$}}
                        \end{picture}}
\psfrag{C}[][]{\begin{picture}(0,0)
            \put(-2,-3){\makebox(0,0)[l]{$C$}}
                        \end{picture}}

\psfrag{At}[][]{\begin{picture}(0,0)
            \put(-5,-5){\makebox(0,0)[l]{$A_{1/2}$}}
                        \end{picture}}

\psfrag{Bt}[][]{\begin{picture}(0,0)
            \put(4,10){\makebox(0,0)[l]{$B_{1/2}$}}
                        \end{picture}}

\psfrag{Ct}[][]{\begin{picture}(0,0)
            \put(-24,12){\makebox(0,0)[l]{$C_{1/2}$}}
                        \end{picture}}

\resizebox{0.482\linewidth}{!}{%
 \includegraphics{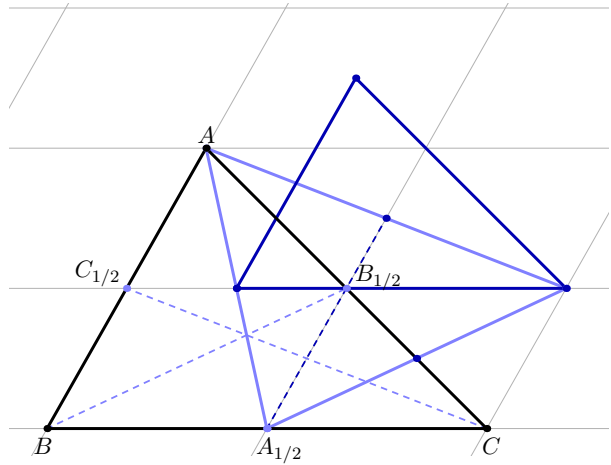}}
    \caption{
A ``proof'' of Property (\ref{Prd})
        }
\label{fsmtmt}

\end{figure}

Proving Property~(\ref{Pra}) is a common exercise. We provide ``proofs'' without words of Properties~(\ref{Prb}), (\ref{Prc}) and~(\ref{Prd}) in Figures~\ref{famt} and~\ref{fsmtmt}.  Different proofs can be browsed at \cite{web1} and \cite{web3}. Note that Property~(\ref{Prb}) fails for other equally important intersecting cevians of a triangle; for example,  as shown in \cite{B03}, one cannot speak about a triangle formed by bisectors or altitudes. As usual, we say that three line segments form a triangle if there exists a triangle whose sides have the same lengths as the given line segments.

\begin{figure}[H]

\setlength{\abovecaptionskip}{5pt}%
\setlength{\belowcaptionskip}{-5pt}%

\psfrag{A}[][]{\begin{picture}(0,0)
            \put(3,4){\makebox(0,0)[l]{$A$}}
                        \end{picture}}

\psfrag{At}[][]{\begin{picture}(0,0)
            \put(-3,-9){\makebox(0,0)[l]{$A_{1/2}$}}
                        \end{picture}}
\psfrag{Ar}[][]{\begin{picture}(0,0)
            \put(-3,-9){\makebox(0,0)[l]{$A_{-1/2}$}}
                        \end{picture}}
\psfrag{As}[][]{\begin{picture}(0,0)
            \put(-3,-9){\makebox(0,0)[l]{$A_{3/2}$}}
                        \end{picture}}

\psfrag{B}[][]{\begin{picture}(0,0)
            \put(-3,-9){\makebox(0,0)[l]{$B$}}
                        \end{picture}}

\psfrag{Bt}[][]{\begin{picture}(0,0)
            \put(-13,-11){\makebox(0,0)[l]{$B_{1/2}$}}
                        \end{picture}}
\psfrag{Br}[][]{\begin{picture}(0,0)
            \put(-19,-9){\makebox(0,0)[l]{$B_{-1/2}$}}
                        \end{picture}}
\psfrag{Bs}[][]{\begin{picture}(0,0)
            \put(0,5){\makebox(0,0)[l]{$B_{3/2}$}}
                        \end{picture}}

\psfrag{C}[][]{\begin{picture}(0,0)
            \put(-6,-9){\makebox(0,0)[l]{$C$}}
                        \end{picture}}

\psfrag{Ct}[][]{\begin{picture}(0,0)
            \put(-2,-10){\makebox(0,0)[l]{$C_{1/2}$}}
                        \end{picture}}
\psfrag{Cr}[][]{\begin{picture}(0,0)
            \put(4,3){\makebox(0,0)[l]{$C_{-1/2}$}}
                        \end{picture}}
\psfrag{Cs}[][]{\begin{picture}(0,0)
            \put(0,-9){\makebox(0,0)[l]{$C_{3/2}$}}
                        \end{picture}}

\psfrag{Ga}[][]{\begin{picture}(0,0)
            \put(-5,-8){\makebox(0,0)[l]{$G_{a}$}}
                        \end{picture}}

\psfrag{Gb}[][]{\begin{picture}(0,0)
            \put(1,6){\makebox(0,0)[l]{$G_{b}$}}
                        \end{picture}}

\psfrag{Gc}[][]{\begin{picture}(0,0)
            \put(-15,5){\makebox(0,0)[l]{$G_{c}$}}
                        \end{picture}}

\resizebox{0.7\linewidth}{!}{%
  \includegraphics{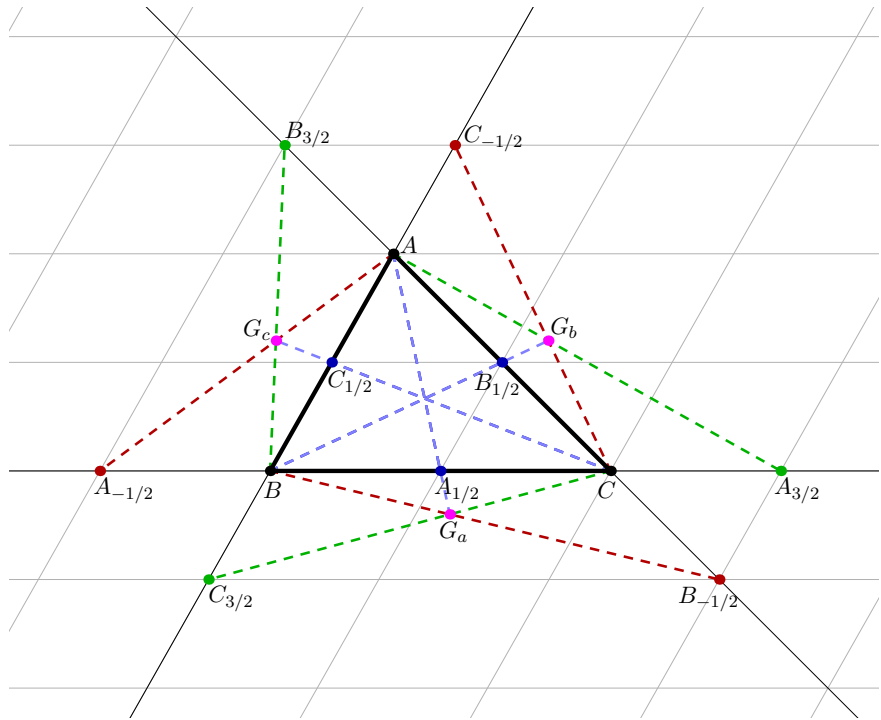}}
    \caption{
The ``median grid''; three medians and six outer medians
        }
\label{fMedGrid}

\end{figure}

A median of a triangle is just a special \emph{cevian}; a cevian is a line segment joining a vertex of a triangle to a point on the opposite side. Are there other triples of cevians from distinct vertices of a triangle that share Properties~(\ref{Pra}), (\ref{Prb}), ~(\ref{Prc}) and~(\ref{Prd})?  Some natural candidates for such cevians are suggested by the ``median grid'' already encountered in Figures~\ref{famt} and~\ref{fsmtmt}.  In Figure~\ref{fMedGrid} we show more of this grid with the cevians that are in some sense closest to the medians.  The labeling of  the points on the line $BC$ in Figure~\ref{fMedGrid} originates from $BC$ being considered as a number line with $0$ at $B$ and $1$ at $C$.  More precisely, for $\rho \in \nR$, the point $A_\rho$ is the point on the line $BC$ which satisfies $\overrightarrow{B\!A}_\rho =\rho \overrightarrow{BC}$.  The points on the lines $AB$ and $BC$ are labeled  similarly.  As indicated in Figure~\ref{fMedGrid}, the cevians in the following three triples
\begin{equation} \label{eqOMT}
\bigl(AA_{-1/2}, CC_{1/2}, B\!\rule{1pt}{0pt}B_{3/2}\bigr), \quad
 \bigl(B\!\rule{1pt}{0pt}B_{-1/2}, AA_{1/2}, CC_{3/2}\bigr), \quad
  \bigl(CC_{-1/2}, B\!\rule{1pt}{0pt}B_{1/2}, AA_{3/2}\bigr),
\end{equation}
are concurrent; with the understanding that three line segments are concurrent if the lines determined by them are concurrent.  That this is really true can be proved as an exercise in vector algebra or by using Ceva's characterization of concurrent cevians: $AA_{\rho}$, $CC_{\sigma}$, $B\!\rule{1pt}{0pt}B_{\tau}$ are concurrent if and only if
\begin{equation}\label{eqCeva}
\rho\sigma\tau - (1-\rho)(1-\sigma)(1-\tau) = 0.
\end{equation}
Since the cevians $AA_{-1/2}$, $AA_{3/2}$, $B\!\rule{1pt}{0pt}B_{-1/2}$, $B\!\rule{1pt}{0pt}B_{3/2}$,  $CC_{-1/2}$, $CC_{3/2}$ play the leading roles in this note and because of their proximity to the medians on the ``median grid'', we call them \emph{outer medians}. Thus, for example, associated to vertex $A$ we have one median, $AA_{1/2}$, and two outer medians, $AA_{-1/2}$ and $AA_{3/2}$, see  Figure~\ref{fMedGrid}.

\begin{figure}[H]

\setlength{\abovecaptionskip}{5pt}%
\setlength{\belowcaptionskip}{-5pt}%

\psfrag{A}[][]{\begin{picture}(0,0)
            \put(3,4){\makebox(0,0)[l]{$A$}}
                        \end{picture}}

\psfrag{At}[][]{\begin{picture}(0,0)
            \put(-3,-9){\makebox(0,0)[l]{$A_{1/2}$}}
                        \end{picture}}
\psfrag{Ar}[][]{\begin{picture}(0,0)
            \put(-3,-9){\makebox(0,0)[l]{$A_{-1/2}$}}
                        \end{picture}}
\psfrag{As}[][]{\begin{picture}(0,0)
            \put(-3,-9){\makebox(0,0)[l]{$A_{3/2}$}}
                        \end{picture}}

\psfrag{B}[][]{\begin{picture}(0,0)
            \put(-12,-5){\makebox(0,0)[l]{$B$}}
                        \end{picture}}

\psfrag{Bt}[][]{\begin{picture}(0,0)
            \put(-13,-11){\makebox(0,0)[l]{$B_{1/2}$}}
                        \end{picture}}
\psfrag{Br}[][]{\begin{picture}(0,0)
            \put(2,3){\makebox(0,0)[l]{$B_{-1/2}$}}
                        \end{picture}}
\psfrag{Bs}[][]{\begin{picture}(0,0)
            \put(3,3){\makebox(0,0)[l]{$B_{3/2}$}}
                        \end{picture}}

\psfrag{C}[][]{\begin{picture}(0,0)
            \put(-6,-7){\makebox(0,0)[l]{$C$}}
                        \end{picture}}

\psfrag{Ct}[][]{\begin{picture}(0,0)
            \put(-2,-10){\makebox(0,0)[l]{$C_{1/2}$}}
                        \end{picture}}
\psfrag{Cr}[][]{\begin{picture}(0,0)
            \put(4,3){\makebox(0,0)[l]{$C_{-1/2}$}}
                        \end{picture}}
\psfrag{Cs}[][]{\begin{picture}(0,0)
            \put(0,-9){\makebox(0,0)[l]{$C_{3/2}$}}
                        \end{picture}}

\resizebox{0.7\linewidth}{!}{%
  \includegraphics{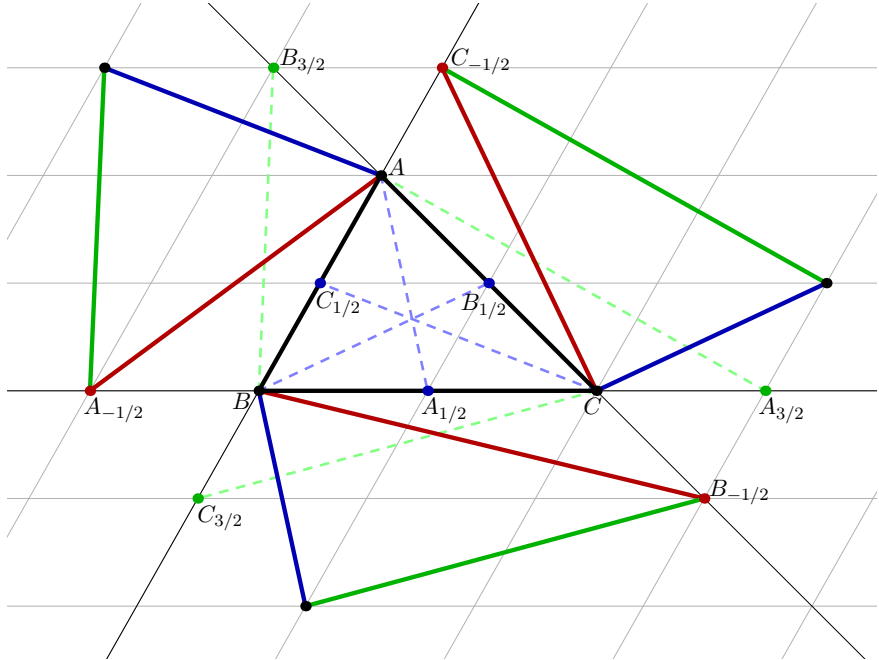}}
    \caption{
Three outer median triangles of $ABC$
        }
\label{fOMT3}

\end{figure}

We find it quite remarkable that all four properties of the medians listed in the opening of this note hold for the three triples displayed in \eqref{eqOMT}, each of which consists of a median and two outer medians originating from distinct vertices.
\begin{enumerate}[{\rm (A)}]
\item  \label{PrA}
The median and two outer medians in each of the triples in \eqref{eqOMT} are concurrent. The concurrency point divides the outer medians in the ratio 2:3.

\item \label{PrB}
The median and two outer medians in the triples in \eqref{eqOMT} form three triangles. We refer to these three triangles as \emph{outer median triangles} of $ABC$, see Figure~\ref{fOMT3}.

\item \label{PrC}
The area of each outer median triangle of $ABC$ is 5/4 of the area of $ABC$.

\item \label{PrD}
For each outer median triangle, one of its outer median triangles is similar to the original triangle $ABC$ with the ratio of similarity 5/4.
\end{enumerate}

\begin{figure}[H]

\setlength{\abovecaptionskip}{5pt}%
\setlength{\belowcaptionskip}{-5pt}%

\psfrag{A}[][]{\begin{picture}(0,0)
            \put(-3,2){\makebox(0,0)[l]{$A$}}
                        \end{picture}}

\psfrag{At}[][]{\begin{picture}(0,0)
            \put(-7,-3){\makebox(0,0)[l]{$A_{1/2}$}}
                        \end{picture}}

\psfrag{B}[][]{\begin{picture}(0,0)
            \put(-4,-3){\makebox(0,0)[l]{$B$}}
                        \end{picture}}

\psfrag{Br}[][]{\begin{picture}(0,0)
            \put(-2,-3){\makebox(0,0)[l]{$B_{-1/2}$}}
                        \end{picture}}

\psfrag{C}[][]{\begin{picture}(0,0)
            \put(-5,-3){\makebox(0,0)[l]{$C$}}
                        \end{picture}}

\psfrag{Cs}[][]{\begin{picture}(0,0)
            \put(-28,-8){\makebox(0,0)[l]{$C_{3/2}$}}
                        \end{picture}}

\resizebox{0.77\linewidth}{!}{%
  \includegraphics{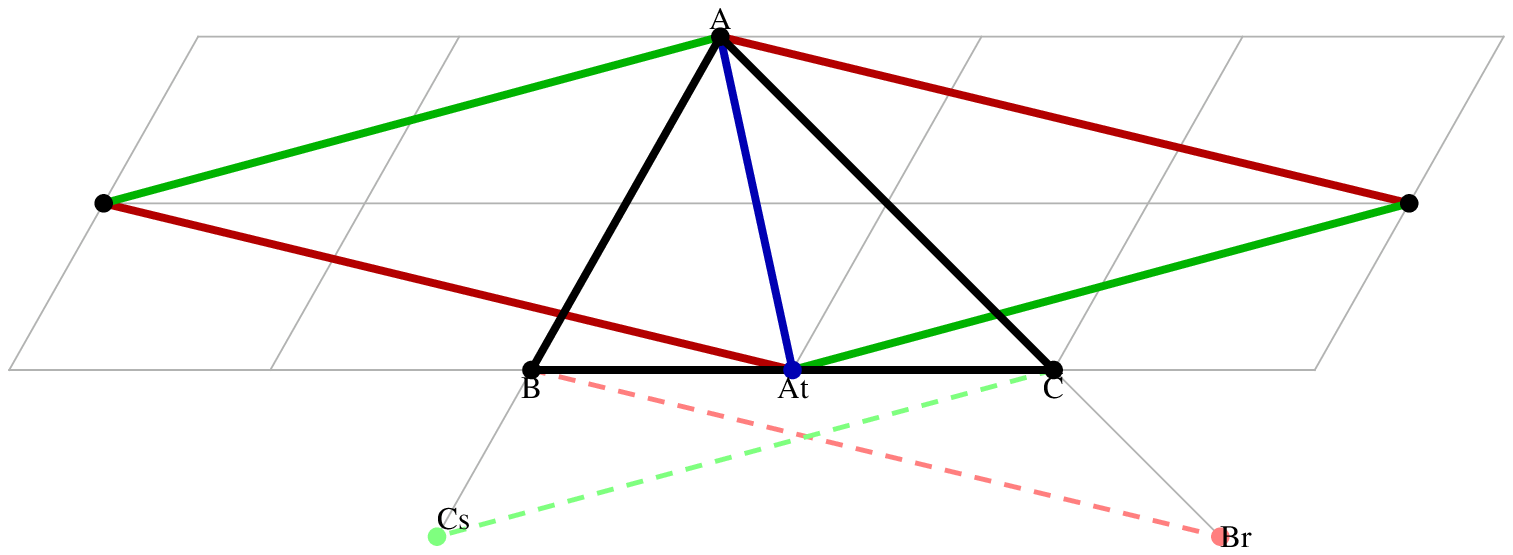}}
    \caption{
A ``proof'' of Properties~(\ref{PrB}) and (\ref{PrC})
        }
\label{faomt}

\end{figure}

\begin{figure}[H]

\setlength{\abovecaptionskip}{5pt}%
\setlength{\belowcaptionskip}{-5pt}%

\psfrag{A}[][]{\begin{picture}(0,0)
            \put(3,1){\makebox(0,0)[l]{$A=Q$}}
                        \end{picture}}

\psfrag{At}[][]{\begin{picture}(0,0)
            \put(-6,-4){\makebox(0,0)[l]{$A_{1/2}\!=\!R$}}
                        \end{picture}}

\psfrag{Qr}[][]{\begin{picture}(0,0)
            \put(3,-4){\makebox(0,0)[l]{$Q_{-1/2}$}}
                        \end{picture}}

\psfrag{Rs}[][]{\begin{picture}(0,0)
            \put(0,12){\makebox(0,0)[l]{$R_{3/2}$}}
                        \end{picture}}

\psfrag{P}[][]{\begin{picture}(0,0)
            \put(-3,-11){\makebox(0,0)[l]{$P$}}
                        \end{picture}}
\psfrag{Pt}[][]{\begin{picture}(0,0)
            \put(-16,-2){\makebox(0,0)[l]{$P_{1/2}$}}
                        \end{picture}}

\psfrag{B}[][]{\begin{picture}(0,0)
            \put(-4,-3){\makebox(0,0)[l]{$B$}}
                        \end{picture}}
\psfrag{Br}[][]{\begin{picture}(0,0)
            \put(-6,-9){\makebox(0,0)[l]{$B_{-1/2}$}}
                        \end{picture}}

\psfrag{C}[][]{\begin{picture}(0,0)
            \put(-5,-3){\makebox(0,0)[l]{$C$}}
                        \end{picture}}

\psfrag{Cs}[][]{\begin{picture}(0,0)
            \put(-6,-15){\makebox(0,0)[l]{$C_{3/2}$}}
                        \end{picture}}

\psfrag{Z}[][]{\begin{picture}(0,0)
            \put(-6,3){\makebox(0,0)[l]{$X$}}
                        \end{picture}}

\resizebox{0.6\linewidth}{!}{%
  \includegraphics{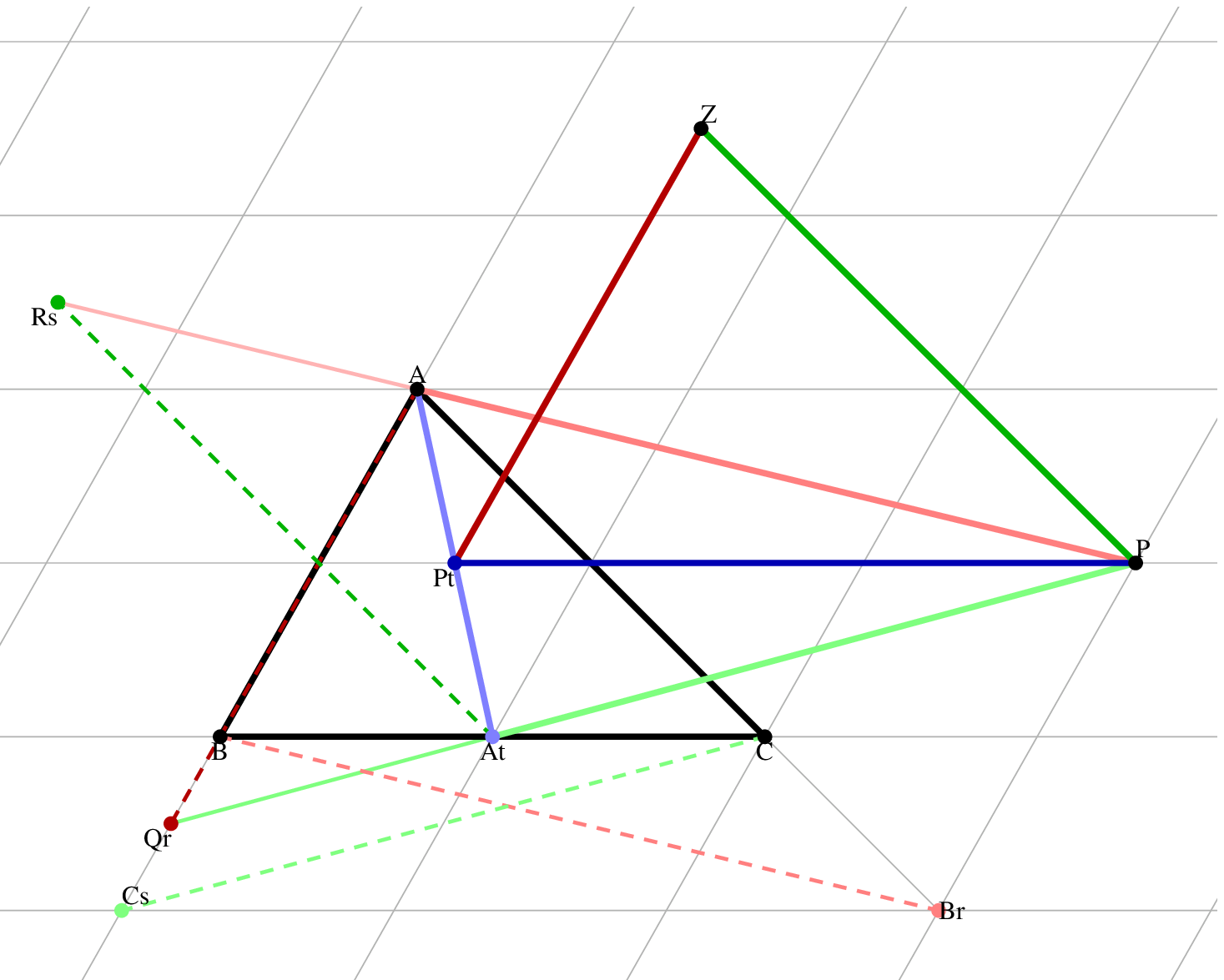}}
    \caption{
A ``proof'' of Property~(\ref{PrD})
        }
\label{fsomt}

\end{figure}

As we already mentioned, Property~(\ref{PrA}) follows from Ceva's theorem, while Figures~\ref{faomt} and ~\ref{fsomt} offer ``proofs'' without words of Properties~(\ref{PrB}), (\ref{PrC}) and~(\ref{PrD}).

We motivated our study of outer medians by their special position on the ``median grid''. However, the above four properties could very well hold for other triples of cevians.  Is it then the case that the median triangle and the three outer median triangles are truly special? Since Property~(\ref{PrA}) is characterized by \eqref{eqCeva}, we explore Property~(\ref{PrB}) next.  We first answer the following question: For which $\rho, \sigma, \tau \in \nR$ does there exist a triangle with the sides which are congruent and parallel to the cevians $AA_\rho, BB_\sigma$ and $CC_\tau$, independent of the triangle $ABC$ in which they are constructed?

Define the vectors $\Vb{a} = \overrightarrow{BC}$, $\Vb{b}  = \overrightarrow{C\!A}$ and
$\Vb{c}  = \overrightarrow{AB}$. We have
\begin{equation} \nonumber  % \label{eqVCevs}
 \overrightarrow{AA}_\rho = \Vb{c} + \rho \, \Vb{a}, \quad
 \overrightarrow{BB}_\sigma = \Vb{a} + \sigma\, \Vb{b} \quad
 \text{and} \quad
 \overrightarrow{CC}_\tau = \Vb{b} + \tau \, \Vb{c}.
\end{equation}
Then, a necessary and sufficient condition for the existence of a triangle with the sides which are congruent and parallel to the line segments $AA_\rho, BB_\sigma$ and $CC_\tau$ is that one of the following four vector equations is satisfied:
\begin{equation} \label{eqVs}
 \overrightarrow{AA}_\rho \hat{\pm} \overrightarrow{B\!\rule{1pt}{0pt}B}_\sigma \mathbf{\pm}
 \overrightarrow{CC}_\tau = \Vb{0}.
\end{equation}
We put a special sign $\hat{\phantom{\cdot}}$ above the first $\pm$ in order to be able to trace this sign in the calculations that follow. Substituting $\Vb{c} = - \Vb{a} - \Vb{b}$ in \eqref{eqVs} we get
\begin{equation} \label{eqVs1}
\bigl(-1 \hat{\pm} 1 + \rho \mp \tau \bigr) \Vb{a}
 + \bigl( -1 \hat{\pm} \sigma \pm 1 \mp \tau \bigr)\Vb{b} = \Vb{0}.
\end{equation}
Using the linear independence of $\Vb{a}$ and $\Vb{b}$ and choosing both $+$ signs in \eqref{eqVs}, it follows from \eqref{eqVs1} that $\rho =  \sigma = \tau$.  Choosing the first sign in \eqref{eqVs} to be $+$ and the second to be $-$ we get that $\rho = -\tau, \sigma = 2-\tau$.  Choosing the first sign in \eqref{eqVs} to be $-$ and the second to be $+$ we get $\rho = 2-\sigma$, $\tau = - \sigma$ and choosing both $-$ signs in \eqref{eqVs} we get $\sigma = - \rho$, $\tau = 2 - \rho$.  Thus, we have identified four sets of parameters $(\rho, \sigma, \tau)$ for which, independent of $ABC$, there exists a triangle, possibly degenerate, with the sides which are congruent and parallel to the cevians $AA_\rho$, $B\!\rule{1pt}{0pt}B_\sigma$ and $CC_\tau$:
\begin{equation} \label{eqGood}
(\xi, \xi, \xi), \quad (2-\xi, \xi, -\xi), \quad (-\xi, 2-\xi, \xi), \quad (\xi, -\xi, 2-\xi),  \quad  \xi \in \nR.
\end{equation}

The only concern here is that the cevians $AA_\rho, B\!\rule{1pt}{0pt}B_\sigma$ and $CC_\tau$ might be parallel. However, the condition for the cevians to be parallel is easily established as follows. Since the vector $\overrightarrow{CC}_\tau$ is nonzero, we look for $\lambda, \mu \in \nR$ such that
\begin{equation} \label{eqli}
\Vb{c} + \rho \, \Vb{a} =  \lambda (\Vb{b} + \tau \, \Vb{c}) \quad  \text{and} \quad
\Vb{a} + \sigma\, \Vb{b} =  \mu (\Vb{b} + \tau \, \Vb{c}).
\end{equation}
Substituting $\Vb{c} = - \Vb{a} - \Vb{b}$ in \eqref{eqli} and using the linear independence of $\Vb{a}$ and $\Vb{b}$, we get from the first equation $\lambda = 1/(\tau-1)$, $\rho = 1/(1-\tau)$ and from the second equation $\mu = -1/\tau$, $\sigma = 1-1/\tau$. Hence, the line segments ${AA}_\rho$, ${B\!\rule{1pt}{0pt}B}_\sigma$ and ${CC}_\tau$ are parallel if and only if
\begin{equation} \label{eqparce}
\rho = \frac{1}{1-\xi}, \quad \sigma = 1- \frac{1}{\xi}, \quad \tau = \xi, \quad \xi \in \nR\setminus\{0,1\}.
\end{equation}

\begin{figure}[H]

\setlength{\abovecaptionskip}{2pt}%
\setlength{\belowcaptionskip}{-10pt}%

\psfrag{rho}[][]{\begin{picture}(0,0)
            \put(-17,0){\makebox(0,0)[l]{$\rho$}}
                        \end{picture}}

\psfrag{sigma}[][]{\begin{picture}(0,0)
            \put(15,2){\makebox(0,0)[l]{$\sigma$}}
                        \end{picture}}

\psfrag{tau}[][]{\begin{picture}(0,0)
            \put(3,7){\makebox(0,0)[l]{$\tau$}}
                        \end{picture}}

\psfrag{E}[][]{\begin{picture}(0,0)
            \put(-5,1){\makebox(0,0)[l]{$\nE$}}
                        \end{picture}}

\psfrag{F}[][]{\begin{picture}(0,0)
            \put(-2,2){\makebox(0,0)[l]{$\nF$}}
                        \end{picture}}

\psfrag{G}[][]{\begin{picture}(0,0)
            \put(-5,1){\makebox(0,0)[l]{$\nG$}}
                        \end{picture}}

\resizebox{0.7\linewidth}{!}{%
  \includegraphics{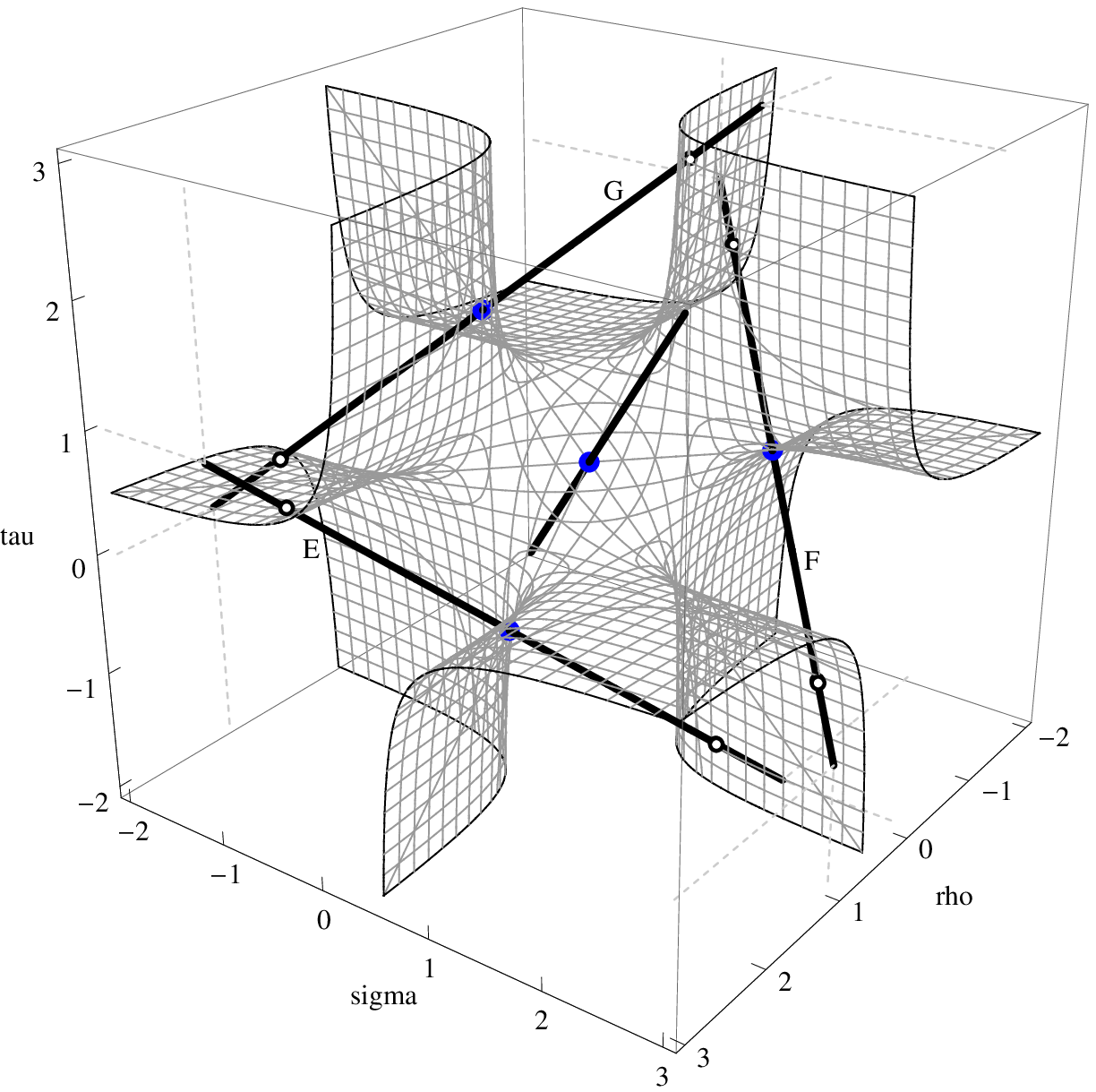}}
    \caption{The sets $\nD$, $\nE$, $\nF$ and $\nG$ and the Ceva surface}
\label{SetA}
\end{figure}

Thus, to avoid degeneracy of triangles with cevian sides corresponding to triples in \eqref{eqGood} such as, for example, $(-\xi, 2-\xi, \xi)$, we must exclude the values of the parameter $\xi$ that solve $-\xi=\tfrac{1}{1-\xi}$. This, in turn, shows that the triples $(\rho, \sigma, \tau)$ for which there exists a non-degenerate triangle with the sides which are congruent and parallel to the cevians $AA_\rho$, $B\!\rule{.5pt}{0pt}B_\sigma$ and $CC_\tau$ must belong to one of the following four sets:
\begin{alignat*}{2}
\nD & = \bigl\{ (\xi, \xi, \xi)  :  \xi \in \nR \bigr\}, \quad &
 \nE & =
 \bigl\{  (2-\xi, \xi, -\xi)  :   \xi \in \nR\!\setminus\!\{-\phi^{-1},\phi\} \bigr\}, \\
 \nF & = \bigl\{ (-\xi, 2-\xi, \xi) :   \xi \in \nR\!\setminus\!\{-\phi^{-1},\phi\} \bigr\},  \quad &
  \nG & = \bigl\{ (\xi, -\xi, 2-\xi)   :   \xi \in \nR\!\setminus\!\{-\phi^{-1},\phi\} \bigr\},
\end{alignat*}
where $\phi=(1+\sqrt{5})/2$ denotes the golden ratio. The diagonal of  the $\rho\sigma\tau$-space provides a geometric representation of the set $\nD$,  while the other three sets are represented by straight lines with two points removed.  These lines are shown in Figure~\ref{SetA} together with the surface whose equation is \eqref{eqCeva}, which we call the {\em Ceva surface}.

\begin{figure}[H]

\setlength{\abovecaptionskip}{5pt}%
\setlength{\belowcaptionskip}{-5pt}%

\psfrag{A}[][]{\begin{picture}(0,0)
            \put(-3,2){\makebox(0,0)[l]{$A$}}
                        \end{picture}}
\psfrag{At}[][]{\begin{picture}(0,0)
            \put(-6,-4){\makebox(0,0)[l]{$A_{\xi}$}}
                        \end{picture}}

\psfrag{B}[][]{\begin{picture}(0,0)
            \put(-5,-2){\makebox(0,0)[l]{$B$}}
                        \end{picture}}

\psfrag{Bt}[][]{\begin{picture}(0,0)
            \put(-3,-2){\makebox(0,0)[l]{$B_{\xi}$}}
                        \end{picture}}

\psfrag{C}[][]{\begin{picture}(0,0)
            \put(-3,-2){\makebox(0,0)[l]{$C$}}
                        \end{picture}}

\psfrag{Ct}[][]{\begin{picture}(0,0)
            \put(-10,-2){\makebox(0,0)[l]{$C_{\xi}$}}
                        \end{picture}}

\resizebox{0.525\linewidth}{!}{%
  \includegraphics{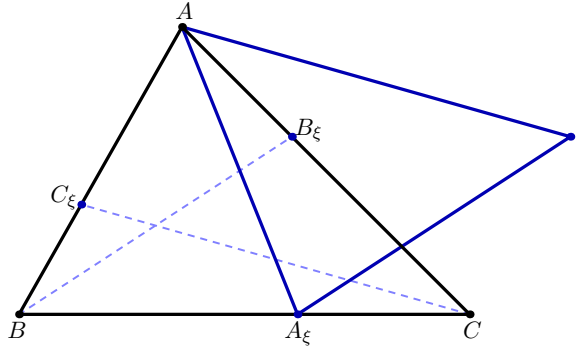}}
    \caption{
The $\xi$-median triangle with $\xi = \phi^{-1}$
        }
\label{fxMT}

\end{figure}

\begin{figure}[H]

\setlength{\abovecaptionskip}{5pt}%
\setlength{\belowcaptionskip}{-5pt}%

\psfrag{A}[][]{\begin{picture}(0,0)
            \put(3,4){\makebox(0,0)[l]{$A$}}
                        \end{picture}}

\psfrag{At}[][]{\begin{picture}(0,0)
            \put(-3,-9){\makebox(0,0)[l]{$A_{\xi}$}}
                        \end{picture}}
\psfrag{Ar}[][]{\begin{picture}(0,0)
            \put(-7,-9){\makebox(0,0)[l]{$A_{-\xi}$}}
                        \end{picture}}
\psfrag{As}[][]{\begin{picture}(0,0)
            \put(-5,-9){\makebox(0,0)[l]{$A_{2-\xi}$}}
                        \end{picture}}

\psfrag{B}[][]{\begin{picture}(0,0)
            \put(-13,-6){\makebox(0,0)[l]{$B$}}
                        \end{picture}}

\psfrag{Bt}[][]{\begin{picture}(0,0)
            \put(-6,-10){\makebox(0,0)[l]{$B_{\xi}$}}
                        \end{picture}}
\psfrag{Br}[][]{\begin{picture}(0,0)
            \put(2,2){\makebox(0,0)[l]{$B_{-\xi}$}}
                        \end{picture}}
\psfrag{Bs}[][]{\begin{picture}(0,0)
            \put(2,4){\makebox(0,0)[l]{$B_{2-\xi}$}}
                        \end{picture}}

\psfrag{C}[][]{\begin{picture}(0,0)
            \put(-6,-8){\makebox(0,0)[l]{$C$}}
                        \end{picture}}

\psfrag{Ct}[][]{\begin{picture}(0,0)
            \put(-1,-9){\makebox(0,0)[l]{$C_{\xi}$}}
                        \end{picture}}
\psfrag{Cr}[][]{\begin{picture}(0,0)
            \put(-22,2){\makebox(0,0)[l]{$C_{-\xi}$}}
                        \end{picture}}
\psfrag{Cs}[][]{\begin{picture}(0,0)
            \put(-25,4){\makebox(0,0)[l]{$C_{2-\xi}$}}
                        \end{picture}}

\resizebox{0.7\linewidth}{!}{%
  \includegraphics{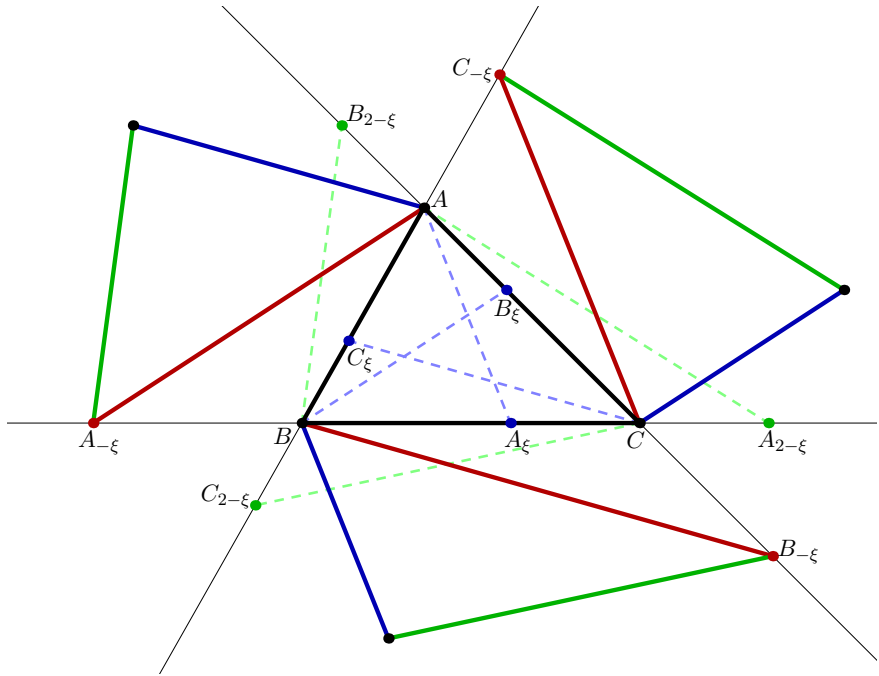}}
    \caption{
Three $\xi$-outer median triangles with $\xi = \phi^{-1}$
        }
\label{fxOMT3}

\end{figure}

Before returning to Figure~\ref{SetA}, we will explore the large family of triples of cevians determined by the triples in the sets $\nD$, $\nE$, $\nF$ and $\nG$.  The cevians associated with these triples are guaranteed to form triangles, that is  they satisfy a property analogous to Property~(\ref{PrB}). The most prominent representatives of triangles originating from the sets $\nD$, $\nE$, $\nF$ and $\nG$ are the median and outer median triangles, which all correspond to the value $\xi = 1/2$. Therefore, for a fixed $\xi$, the triangle associated with the triple $(\xi,\xi,\xi)$ in $\nD$ we call {\em$\xi$-median triangle} and the triangles associated with the corresponding triples in $\nE$, $\nF$ and $\nG$ we call {\em$\xi$-outer median triangles}. In Figures~\ref{fxMT} and~\ref{fxOMT3} we illustrate these triangles with $\xi = 1/\phi$,  the reciprocal of the golden ratio.

Next, we explore whether the $\xi$-median  and  $\xi$-outer median triangles  have properties analogous to~(\ref{PrC}) and~(\ref{PrD}). Two classical formulas seem to be custom made for this task. First, Heron's formula which gives the square of the area of a triangle in terms of squares of its sides $a, b, c$:
\begin{equation} \nonumber
\frac{1}{16} \Bigl( 2\bigl(a^2 b^2 + b^2 c^2 + c^2 a^2 \bigr) -  \bigl(a^4+b^4+c^4\bigr) \Bigr).
\end{equation}
Second, Stewart's theorem which gives the square of the length of a cevian in terms of the squares of the sides of $ABC$. Specifically, in matrix form, these equations are:
\begin{equation} \label{eqM}
\left[\!\!\begin{array}{c}
\bigl(CC_\tau\bigr)^2 \\[6pt]
 \bigl(B\!\rule{1pt}{0pt}B_\sigma\bigr)^2 \\[6pt]
  \bigl(AA_\rho\bigr)^2
\end{array}\!\!\right] = \left[\!
\begin{array}{ccc}
 \tau & 1-\tau  & \tau (\tau -1) \\[6pt]
  1-\sigma  &  \sigma (\sigma -1)  & \sigma \\[6pt]
  \rho (\rho - 1)  & \rho  &  1-\rho
\end{array}\!\right] \ \left[\!\begin{array}{c}
a^2 \\[6pt]
 b^2 \\[6pt]
  c^2
\end{array}\!\!\right].
\end{equation}
The idea to write this operation in the above matrix form is due to Griffiths  \cite{Griffiths}.  We denote the $3\times 3$ matrix in \eqref{eqM} by $M(\rho,\sigma,\tau)$.

Now it is clear how to proceed to verify the property analogous to~(\ref{PrC}): use triples from the sets $\nD$, $\nE$, $\nF$ and $\nG$ to get expressions for the squares of the corresponding cevians, substitute these expressions in Heron's formula and simplify. However, this involves simplifying an expression with 36 additive terms; quite a laborious task for a human but a perfect challenge for a Computer Algebra System like {\em Mathematica}.  We first define Heron's formula as a {\em Mathematica} function, we call it  {\tt HeronS}, operating on the {\em triples} of squares of the sides of a triangle and producing the square of the area:

\resizebox{!}{!}{%
\includegraphics{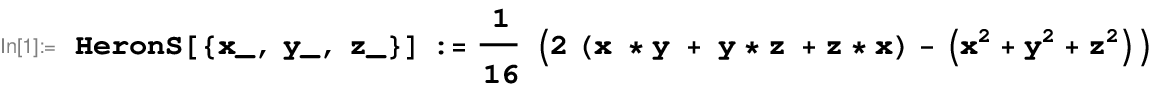}}

\noindent Next, we define in {\em Mathematica} the matrix function {\tt M} as in \eqref{eqM}:

\rule{0pt}{38pt}\resizebox{1\width}{!}{%
\includegraphics{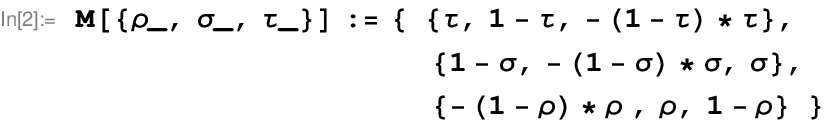}}

\vspace*{-4pt}

To verify the property analogous to (\ref{PrC}) for $\xi$-median triangles, we put the newly defined functions in action by calculating the ratio between the squares of the area of the $\xi$-median triangle and the original triangle. {\em Mathematica}'s answer is instantaneous:

\rule{0pt}{37pt}\resizebox{!}{!}{%
\includegraphics{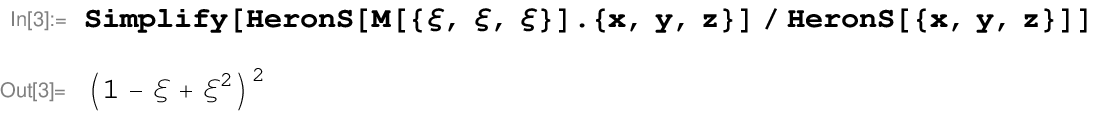}}

\vspace*{-4pt}
\noindent This ``proves'' that the ratio of the areas depends only on $\xi$ and it is $1-\xi + \xi^2$. Further, for one of the $\xi$-outer median triangles we have

\rule{0pt}{36pt}\resizebox{!}{!}{%
\includegraphics{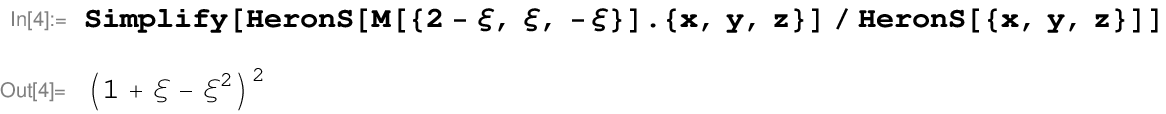}}

\vspace*{-4pt}
\noindent ``proving'' that the area of the triangle formed by the cevians $AA_{2-\xi}$, $B\!\rule{1pt}{0pt}B_{\xi}$, $CC_{-\xi}$ is $|1+\xi - \xi^2|$ of the area of the original triangle $ABC$. The other two $\xi$-outer median triangles yield the same ratio.  In summary, {\em Mathematica} has confirmed that the $\xi$-median and the three $\xi$-outer median triangles all have the property analogous to (\ref{PrC}).

The verification of the property analogous to (\ref{PrD}) is simpler. For a $\xi$-median triangle, following \cite{Griffiths}, we just need to calculate the square of the matrix $M(\xi,\xi,\xi)$, which turns out to be $\bigl(1-\xi+\xi^2\bigl)^2 I$.  This confirms that the $\xi$-median triangle of the $\xi$-median triangle is similar to the original triangle with the ratio of similarity $1-\xi+\xi^2$.

Similarly, for a $\xi$-outer median triangle corresponding to a triple in $\nE$, we calculate the square of the matrix $M(2-\xi,\xi,-\xi)$ which turns out to be $\bigl(1+\xi-\xi^2\bigl)^2 I$; thus confirming that one of the $\xi$-outer median triangles of this $\xi$-outer median triangle is similar to the original triangle with the ratio of similarity $|1+\xi-\xi^2|$. In contrast, for a $\xi$-outer median triangle corresponding to a triple in $\nF$, to get a triangle similar to the original triangle we need to calculate its $\xi$-outer median triangle corresponding to a triple in $\nG$. This amounts to multiplying the matrices
\[
M(\xi,-\xi,2-\xi)M(-\xi,2-\xi,\xi) = \bigl(1+\xi-\xi^2\bigl)^2 I.
\]
Likewise, for a $\xi$-outer median triangle corresponding to a triple in $\nG$, we calculate its $\xi$-outer median triangle corresponding to a triple in $\nF$ and obtain the same result.

All these calculations indicate that, after all, the median and outer median triangles are facing a stiff competition from their $\xi$-triangles  generalizations.  However, property (\ref{PrA}) comes to the rescue of the median and outer median triangles at this point. We want the triples of cevians corresponding to the triples in $\nD$, $\nE$, $\nF$ and $\nG$ to be concurrent as well.  So, which of these triples satisfy Ceva's condition \eqref{eqCeva}?  Or, geometrically, what is the intersection of the lines and the Ceva surface in Figure~\ref{SetA}?  First, we substitute $\rho = \sigma = \tau = \xi$ in \eqref{eqCeva}, which yields $\xi^3-(1-\xi)^3=0$, whose only real solution is $\xi=1/2$.  The corresponding cevians are the medians.  To intersect $\nE$ with the Ceva surface we substitute $(2-\xi, \xi, -\xi)$ in \eqref{eqCeva}, obtaining $-\xi^2 (2-\xi) - (1+\xi)(\xi-1)(1-\xi) = 0$, which is equivalent to $(\xi-\phi)(\xi+\phi^{-1})(2\xi-1) = 0$. Since $\xi\not\in\{\phi, -\phi^{-1}\}$, the only solution is $\xi=1/2$, yielding the ``outer median triple'' $(3/2, 1/2,-1/2)$. Intersecting $\nF$ with the Ceva surface gives  $(-1/2, 3/2, 1/2)$ and intersecting $\nG$ with the Ceva surface results in $(1/2, -1/2, 3/2)$.  Consequently, the only triples in $\nD$, $\nE$, $\nF$ and $\nG$ which correspond to concurrent cevians are the ``median triple'' and the three ``outer median triples''.

There is only a slight weakness in our argument above. In identifying the sets  $\nD$, $\nE$, $\nF$ and $\nG$ we assumed that the triangles formed by the corresponding cevians have sides that are {\em parallel} to the cevians themselves.  In \cite{BC13} we proved that the only cevians $AA_\rho, BB_\sigma, CC_\tau$ that form triangles and with $(\rho, \sigma, \tau)$ not included in the sets $\nD$, $\nE$, $\nF$ and $\nG$ are parallel cevians, that is the cevians ${AA}_\rho$, ${B\!\rule{1pt}{0pt}B}_\sigma$ and ${CC}_\tau$, where $\rho$, $\sigma$, $\tau$ satisfy \eqref{eqparce} with the additional restriction
\begin{equation*}
\xi \in  \bigl(  -\phi, -\phi^{-1} \bigr) \! \cup \!
  \bigl( \phi^{-2}, \phi^{-1} \bigr)\! \cup\!
  \bigl( \phi, \phi^2 \bigr).
\end{equation*}
As it turns out, the properties analogous to (\ref{PrC}) and (\ref{PrD}) do not hold for triangles formed by such cevians.  In conclusion, indeed, along with the medians and the median triangle, the outer medians and their outer median triangles are unique in satisfying  all four properties  analogous to those from the beginning of our note.

\end{document}